\documentclass[12pt,a4paper]{article}

\usepackage{stmaryrd}

\usepackage{hyperref}

\newcommand{\pageformat}[6]{\setlength{\hoffset}{-1in}
                  \setlength{\voffset}{-1in}
                  \addtolength{\hoffset}{#5}
                            \addtolength{\voffset}{#6}
                            \setlength{\oddsidemargin}{#1}
                            \setlength{\evensidemargin}{#2}
                            \setlength{\textwidth}{\paperwidth}
                  \addtolength{\textwidth}{-\oddsidemargin}
                  \addtolength{\textwidth}{-\evensidemargin}
                  \addtolength{\textwidth}{-\marginparsep}
                  \addtolength{\textwidth}{-\marginparwidth}
                            \setlength{\topmargin}{#3}
                            \setlength{\textheight}{\paperheight}
                  \addtolength{\textheight}{-\topmargin}
                  \addtolength{\textheight}{-\headheight}
                  \addtolength{\textheight}{-\headsep}
                  \addtolength{\textheight}{-\footskip}
                  \addtolength{\textheight}{-#4}}
\pageformat{2.5cm}{2.5cm}{25mm}{25mm}{1pt}{0pt}

\usepackage{ifthen}
\newboolean{article}
    \setboolean{article}{true}
\newboolean{report}
\newboolean{book}
\newboolean{letter}
\newboolean{german}
\newboolean{italian}
\newboolean{nobaselinestretch}
\newboolean{nosectionappendix}
\newboolean{oldtoc}
\newboolean{nosectionequn}
\newboolean{notheorem}

\ifthenelse{\boolean{german}}{
    \usepackage{german}}{}

\usepackage[latin1]{inputenc}

\makeatletter
\ifthenelse{\boolean{article}\and\not\boolean{oldtoc}}{
    \renewcommand*\l@section[2]{%
        \ifnum \c@tocdepth >\z@
            \addpenalty\@secpenalty
            \addvspace{0.0em \@plus\p@}%
            \setlength\@tempdima{1.5em}%
            \begingroup
                \parindent \z@ \rightskip \@pnumwidth
                \parfillskip -\@pnumwidth
                \leavevmode 
                \advance\leftskip\@tempdima
                \hskip -\leftskip
                #1\nobreak\dotfill \nobreak\hb@xt@\@pnumwidth{\hss #2}\par
            \endgroup
        \fi}}{}
\makeatother

\usepackage{amsmath}
\usepackage{amssymb}
\usepackage[mathscr]{eucal}

\ifthenelse{\boolean{notheorem}}{}{
    \usepackage{theorem}}

\ifthenelse{\boolean{nobaselinestretch}}{}{
    \renewcommand{\baselinestretch}{1.25}}

\newenvironment{env}[2]{\begin{#1}#2\end{#1}}{}
    \newcommand{\beq}[1]{\begin{env}{equation}{#1}}
    \newcommand{\beqn}[1]{\begin{env}{equation*}{#1}}
    \newcommand{\bal}[1]{\begin{env}{align}{#1}}
    \newcommand{\baln}[1]{\begin{env}{align*}{#1}}
    \newcommand{\bga}[1]{\begin{env}{gather}{#1}}
    \newcommand{\bgan}[1]{\begin{env}{gather*}{#1}}
    \newcommand{\bflal}[1]{\begin{env}{flalign}{#1}}
    \newcommand{\bflaln}[1]{\begin{env}{flalign*}{#1}}
    \newcommand{\bmu}[1]{\begin{env}{multline}{#1}}
    \newcommand{\bmun}[1]{\begin{env}{multline*}{#1}}
    \newcommand{\bsp}[1]{\begin{env}{split}{#1}}

    \newcommand{\eeq}{\end{env}}
    \newcommand{\eeqn}{\end{env}}
    \newcommand{\eal}{\end{env}}
    \newcommand{\ealn}{\end{env}}
    \newcommand{\ega}{\end{env}}
    \newcommand{\egan}{\end{env}}
    \newcommand{\eflal}{\end{env}}
    \newcommand{\eflaln}{\end{env}}
    \newcommand{\emu}{\end{env}}
    \newcommand{\emun}{\end{env}}
    \newcommand{\esp}{\end{env}}

\newcommand{\lf}{\vspace{2ex}}

\renewcommand{\bf}[1]{\textbf{#1}}
\renewcommand{\it}[1]{\textit{#1}}

\renewcommand{\sf}[1]{\textsf{#1}}

\renewcommand{\tt}[1]{\texttt{#1}}
\newcommand{\hl}[1]{\bf{\it{#1}}}
\newcommand{\mrm}[1]{\mathrm{#1}}
\newcommand{\mbf}[1]{\mathbf{#1}}
\newcommand{\msf}[1]{\text{\small $\sf{#1}$}}

\newcommand{\cmc}[1]{\mathcal{#1}}
\newcommand{\eus}[1]{\mathscr{#1}}

\newcommand{\bb}[1]{\mathbb{#1}}

\newcommand{\mtiny}[1]{{\setlength{\arraycolsep}{.3ex}\text{\tiny$#1$}}}
\newcommand{\nbd}[1]{$#1$\nobreakdash--}
\newcommand{\ol}[1]{\overline{#1}}

\newcommand{\abs}[1]{\left\lvert#1\right\rvert}
\newcommand{\norm}[1]{\left\lVert#1\right\rVert}

\newcommand{\AB}[1]{\langle#1\rangle}
\newcommand{\bAB}[1]{\bigl\langle#1\bigr\rangle}

\newcommand{\CB}[1]{\{#1\}}
\newcommand{\bCB}[1]{\bigl\{#1\bigr\}}
\newcommand{\BCB}[1]{\Bigl\{#1\Bigr\}}
\newcommand{\SB}[1]{[#1]}

\newcommand{\Matrix}[1]{\begin{pmatrix}#1\end{pmatrix}}

\newcommand{\tMatrix}[1]{\mtiny{\Matrix{#1}}}

\newcommand{\rtMatrix}[1]{\raisebox{.3ex}{\tMatrix{#1}}}

\newcommand{\set}[2][]{
    \ifthenelse{\equal{#1}{}}{
        \CB{#2}}{
        \CB{#1~|~#2}}}
\newcommand{\bset}[2][]{
    \ifthenelse{\equal{#1}{}}{
        \bCB{#2}}{
        \bCB{#1~|~#2}}}
\newcommand{\Bset}[2][]{
    \ifthenelse{\equal{#1}{}}{
        \BCB{#2}}{
        \BCB{#1~\big|~#2}}}
\newcommand{\zero}{\CB{0}}

\DeclareMathOperator{\id}{\normalfont\msf{id}}

\renewcommand{\ker}{\operatorname{\msf{ker}}}

\newcommand{\N}{\bb{N}}

\newcommand{\cB}{\cmc{B}}

\newcommand{\sB}{\eus{B}}

\newcommand{\U}{\mbf{1}}
\newcommand{\F}{{\mrm{F}}}

\ifthenelse{\boolean{nosectionequn}}{}{
    \numberwithin{equation}{section}
    }

\ifthenelse{\boolean{article}\or\boolean{letter}\or\boolean{nosectionequn}}{
    \setboolean{nosectionappendix}{true}}{}
\ifthenelse{\boolean{nosectionappendix}}{}{
    \renewcommand{\appendix}{
        \chapter*{\appendixname}
        \addcontentsline{toc}{chapter}{\appendixname}
        \renewcommand{\thesection}{\Alph{section}}
        \setcounter{section}{0}}}
   
\ifthenelse{\boolean{report}\or\boolean{book}}{
    }{}

\ifthenelse{\boolean{notheorem}}{}{
        \newcommand{\mnname}{Mathematical note.}
        \newcommand{\enname}{End of the note.}
        \newcommand{\definame}{Definition.}
        \newcommand{\propname}{Proposition.}
        \newcommand{\lemname}{Lemma.}
        \newcommand{\exname}{Example.}
        \newcommand{\exername}{Exercise.}
        \newcommand{\remname}{Remark.}
        \newcommand{\obname}{Observation.}
        \newcommand{\thmname}{Theorem.}
        \newcommand{\corname}{Corollary.}
        \newcommand{\proofname}{Proof.}
    \ifthenelse{\boolean{german}}{
        \renewcommand{\mnname}{Mathematische Notiz.}
        \renewcommand{\enname}{Ende der Notiz.}
        \renewcommand{\exname}{Beispiel.}
        \renewcommand{\exername}{�bung.}
        \renewcommand{\remname}{Bemerkung.}
        \renewcommand{\obname}{Beobachtung.}
        \renewcommand{\thmname}{Satz.}
        \renewcommand{\corname}{Korollar.}
        \renewcommand{\proofname}{Beweis.}}{}
    \ifthenelse{\boolean{italian}}{
        \renewcommand{\mnname}{Nota matematica.}
        \renewcommand{\enname}{Fina della nota.}
        \renewcommand{\definame}{Definizione.}
        \renewcommand{\propname}{Proposizione.}
        \renewcommand{\exname}{Esempio.}
        \renewcommand{\exername}{Esercizio.}
        \renewcommand{\remname}{Nota.}
        \renewcommand{\obname}{Osservazione.}
        \renewcommand{\thmname}{Teorema.}
        \renewcommand{\corname}{Corollario.}
        \renewcommand{\proofname}{Dimostrazione.}

       \renewcommand{\appendixname}{Appendice}

       }{}
    \theoremheaderfont{\normalfont\bfseries}
    \theoremstyle{change}
        \theorembodyfont{\rmfamily}
            \newtheorem{emp}{}[section]
                \newcommand{\bemp}[1][]{
                    \begin{emp}\hskip-\labelsep\bf{#1}\hskip\labelsep}
                \newcommand{\eemp}{\end{emp}}
\newtheorem{itemp}[emp]{}
                \newcommand{\bitemp}[1][]{
                    \begin{itemp}\hskip-\labelsep\bf{#1}\hskip\labelsep\normalfont\itshape}
                \newcommand{\eitemp}{\end{itemp}}
            \newtheorem{mn}[emp]{\mnname}
                \newcommand{\bnm}{\begin{mn}~\begin{quotation}\renewcommand{\baselinestretch}{1}\small\noindent\ignorespaces}
                \newcommand{\enm}{\end{quotation}\hfill\bf{\enname}\end{mn}}
            \newtheorem{ex}[emp]{\exname}
                \newcommand{\bex}{\begin{ex}}
                \newcommand{\eex}{\end{ex}}
            \newtheorem{exer}[emp]{\exername}
                \newcommand{\bexer}{\begin{exer}}
                \newcommand{\eexer}{\end{exer}}
            \newtheorem{defi}[emp]{\definame}
                \newcommand{\bdefi}{\begin{defi}}
                \newcommand{\edefi}{\end{defi}}
            \newtheorem{rem}[emp]{\remname}
                \newcommand{\brem}{\begin{rem}}
                \newcommand{\erem}{\end{rem}}
            \newtheorem{ob}[emp]{\obname}
                \newcommand{\bob}{\begin{ob}}
                \newcommand{\eob}{\end{ob}}
        \theorembodyfont{\normalfont\itshape}
            \newtheorem{thm}[emp]{\thmname}
                \newcommand{\bthm}{\begin{thm}}
                \newcommand{\ethm}{\end{thm}}
            \newtheorem{prop}[emp]{\propname}
                \newcommand{\bprop}{\begin{prop}}
                \newcommand{\eprop}{\end{prop}}
            \newtheorem{cor}[emp]{\corname}
                \newcommand{\bcor}{\begin{cor}}
                \newcommand{\ecor}{\end{cor}}
            \newtheorem{lem}[emp]{\lemname}
                \newcommand{\blem}{\begin{lem}}
                \newcommand{\elem}{\end{lem}}
\newenvironment{empn}[1]{\lf\noindent\bf{#1}\ignorespaces\hskip\labelsep}{\lf}
		\newcommand{\bempn}[1]{\begin{empn}{#1}}
		\newcommand{\eempn}{\end{empn}}
		\newcommand{\bitempn}[1]{\begin{empn}{#1}\normalfont\itshape}
		\newcommand{\eitempn}{\end{empn}}
                \newcommand{\bnmn}{\begin{empn}{\mnname}~\begin{quotation}\renewcommand{\baselinestretch}{1}\small\noindent\ignorespaces}
                \newcommand{\enmn}{\end{quotation}\hfill\bf{\enname}\end{empn}}
		\newcommand{\bexn}{\begin{empn}{\exname}}
		\newcommand{\eexn}{\end{empn}}
		\newcommand{\bexern}{\begin{empn}{\exername}}
		\newcommand{\eexern}{\end{empn}}
		\newcommand{\bdefin}{\begin{empn}{\definame}}
		\newcommand{\edefin}{\end{empn}}
		\newcommand{\bremn}{\begin{empn}{\remname}}
		\newcommand{\eremn}{\end{empn}}
		\newcommand{\bobn}{\begin{empn}{\obname}}
		\newcommand{\eobn}{\end{empn}}
		\newcommand{\bthmn}{\bitempn{\thmname}}
		\newcommand{\ethmn}{\eitempn}

\newcommand{\qedsymbol}{~\rule[-0.35mm]{2mm}{2mm}}
    \newcounter{proof}[emp]
    \newenvironment{Proof}[1]{
        \vspace{1ex}
        \renewcommand{\item}[1][\stepcounter{proof}(\roman{proof})]%
            {##1\hskip\labelsep}
        \noindent\textsc{#1\hskip\labelsep}}{
        \nolinebreak\qedsymbol}
    \newcommand{\proof}[1][\proofname]{
        \begin{Proof}{#1}\ignorespaces}
    \newcommand{\qed}{\end{Proof}}
    \newcommand{\noqed}{
        \renewcommand{\qedsymbol}{}
        \end{Proof}}}
    \ifthenelse{\boolean{italian}}{
        \renewcommand{\proofname}{Dimostrazione.}}{}

\usepackage[varg]{txfonts}
\usepackage{bbm}

\usepackage[matrix,arrow,curve]{xy}

\begin{document}





\bibliographystyle{amsalpha}
\newcommand{\ins}[1]{}

\title{Polar Decomposition for Non-Adjointable Maps}
\author{}
\author{Michael Skeide}
\date{}
\date{December 2025}

{
\renewcommand{\baselinestretch}{1}
\maketitle


\begin{abstract}
\noindent
We give conditions when not necessarily adjointable operators between Hilbert modules allow for a polar decomposition involving not necessarily adjointable partial isometries. While the latter have been introduced and discussed by Shalit and Skeide \cite{ShaSk23,Ske25}, here we are led, as a basic new ingredient, to the notion of not necessarily adjointable operators $a$ that admit a modulus $\abs{a}$, so-called \it{modular} operators.
\end{abstract}



}






\noindent
Many presentations of Hilbert modules only talk about maps on or between them that have adjoints, and many theorems about such maps explicitly assume they are adjointable. For instance, \it{polar decomposition} is usually defined/sought only for adjointable maps $a$ from a Hilbert \nbd{\cB}mod\-ule $E$ to a Hilbert \nbd{\cB}module $F$ (so that $\abs{a}:=\sqrt{a^*a}\in\sB^a(E)$ exists) and where also the \it{partial isometry} $v\colon E\rightarrow F$ such that $a=v\abs{a}$ and $\ker v=\ker a$ is assumed adjointable explicitly. This polar decomposition is of little help -- simply, because it hardly ever exists: It exists if and only if the closed ranges of $a$ and of $a^*$ are complemented.%
\footnote{
Already here, it is not easy to find a reference. Wegge-Ohlson \cite[Proposition 15.3.7]{WeO93} gives several equivalent criteria but only for $F=E$. Lance \cite{Lan95}, who evidently does like polar decomposition for \nbd{C^*}module operators as much as we do, only mentions it in the end of his Chapter 3 without even formulating the criterion for existence as a proposition. We would like to say it like that: Whoever really needs polar decompositions is better off with von Neumann (or \nbd{W^*}) modules Skeide \cite{Ske00b,Ske06b} (Paschke \cite{Pas73} ans Baillet, Denizeau, and Havet \cite{BDH88}), where polar decomposition alway exists. (And it exists not only for operators in $\sB^a(E)$, because the latter is a von Neumann algebra, but the elements $x\in E$ themselves possess a polar decomposition $x=v\abs{x}$ with $v\in E$.)
}

For an adjointable $v$, the condition to be a partial isometry is, of course, $vv^*v=v$ -- or any equivalent condition, but all of these invoking explicitly existence of $v^*$. But it is easy to find an \hl{isometry} (that is, an inner product preserving map $E\rightarrow F$, automatically \nbd{\cB}linear and contractive) that is not adjointable. In fact, an isometry $v\colon E\rightarrow F$ is adjointable if and only if its range $vE$ is complemented in $F$, that is, if there exists a projection $p\in\sB^a(F)$ onto $vE$; see, for instance, Skeide \cite[Proposition 1.5.13]{Ske01}. (It is a recommended exercise to find a proof before looking into that of \cite[Proposition 1.5.13]{Ske01}.) So, for any non-complemented submodule $E$ of $F$ the canonical embedding $E\rightarrow F$ is a non-adjointable isometry. We think a definition of \it{partial isometry} in which not all isometries are partial isometries, too, is not satisfactory.
 
 In Skeide \cite{Ske25} we started examining partial isometries defined following the proposal with Shalit in \cite[Remark 5.10]{ShaSk23}, based on the observation that we first need to define also \it{coisometry} without explicit reference to adjoints. Recall that a \hl{unitary} is an inner product preserving surjection and that the unitaries are the isomorphisms in a category of Hilbert modules over a fixed \nbd{C^*}algebra. Any unitary $u$ is adjointable with $u^*=u^{-1}$.
 
 \bempn{Definition \cite{ShaSk23}.~}
 Let $E$ and $F$ be Hilbert \nbd{\cB}modules.
 
An \hl{isometry} is a contractive \nbd{\cB}module map $E\rightarrow F$ that corestricts to a unitary $E\rightarrow F'$ for some submodule $F'$ of $F$; a \hl{coisometry} is a contractive \nbd{\cB}module map $v\colon E\rightarrow F$ that restricts to a unitary $E'\rightarrow F$ for some submodule $E'$ of $E$.

A \hl{partial isometry} is a map $E\rightarrow F$ that corestricts to a coisometry $E\rightarrow F'$ for some submodule $F'$ of $F$.
\eempn

Taking into account that the submodule $F'$ to show that $v\colon E\rightarrow F$ is an isometry is bound to be $F'=vE$, this definition of \it{isometry} is, obviously, equivalent to the one given earlier. In \cite[Remark 5.10]{ShaSk23}, we substituted it with an apparently stronger definition only to make the definition of \it{isometry} and of \it{coisometry} a proper \it{co-pair} of definitions (motivating, thus, our definition of coisometry). Based on \cite[Lemma B.1]{ShaSk23}, which asserts that a \bf{contractive} \nbd{\cB}linear idempotent $p$ onto a closed submodule of $E$ is a projection (that is, in particular, $p$ is self-adjoint, thus, adjointable), one shows that every coisometry is adjointable (\cite[Corollary B.2]{ShaSk23}). With this, it follows that our definition is equivalent to what one would define assuming adjointability, namely, a map is a coisometry if it is the adjoint of an adjointable isometry; from this it follows that also the submodule $E'$ to show that $v$ is a coisometry is unique and bound to be $(\ker v)^\perp$.

Now for a partial isometry $v\colon E\rightarrow F$ from the properties of coisometries it survives that $v$ still has a complemented kernel so that we may define the \hl{initial projection} $\pi_v$ of $v$ as the projection onto $(\ker v)^\perp$. (It coincides with $w^*w$ where $w\colon E\rightarrow vE$ is the coisometry obtained by corestricting $v$.) Form the properties of isometries it survives that $v$ is adjointable if and only $vE$ is complemented in $F$. Of course, if $v$ is adjointable, then it is a partial isometry in our sense if and only if it is a partial isometry in the usual, more restrictive, sense. For all this, see the discussion following \cite[Definition 3.2]{Ske25}.

So, knowing now what a partial isometry $v\in\sB^r(E,F)$ is, what else does it take to get a decomposition of $a\in\sB^r(E,F)$ as $a=v\abs{a}$? We know what $\abs{a}$ is when $a$ is adjointable, namely, $\abs{a}=\sqrt{a^*a}$. So, what it actually takes to write down $\abs{a}$ is a positive operator $b$ substituting $a^*a$ when $a^*$ does not exists, because once we have that positive operator, also its positive square root is granted in the \nbd{C^*}algebra $\sB^a(E)$. But for getting a good candidate for $a^*a$ is not indispensable that $a$ be adjointable. In fact, it is enough to find $b\colon E\rightarrow E$ such that
\beqn{ \tag{$*$} \label{*}
\AB{x,by}
~=~
\AB{ax,ay}
}\eeqn
for all $x,y\in E$.%
\footnote{
In \cite{Ske25}, we used this trick to transform equations concerning $a^*a\in\sB^a(E)$ into something that can be checked also if $a$ is not adjointable, frequently. For instance, in \cite[Footnote l]{Ske25} we showed that a projection is a map $p\colon E\rightarrow E$ (\it{a priori} neither linear nor adjointable) that satisfies $\AB{x,py}=\AB{px,py}$ for all $x,y\in E$.
}
Such $b$ is, clearly, \hl{positive} (that is, $\AB{x,bx}\ge0$ for all $x\in E$), self-adjoint, thus, \nbd{\cB}linear and closeable, and, since $E$ is complete, bounded. It may exist, even if $a$ is not adjointable.

\bexn
Suppose $E$ is a non-complemented submodule of $F$ and $a\colon E\rightarrow F$ the canonical injection, then $a$ is not adjointable, but $b=\id_E$ satisfies \eqref{*}. Embedding $a$ and $b$ in the appropriate corners of $\sB^r\rtMatrix{E\\F}$, we also get a pair $a$ and $b$ of operators \bf{on} the module $\rtMatrix{E\\F}=E\oplus F$.
\eexn

Of course, if such $b$ exists, then it is unique.

\bdefin
A map $a\in\sB^r(E,F)$ is \hl{modular} if there exists $b\colon E\rightarrow E$ satisfying \eqref{*}. If $a$ is modular the we define its \hl{modulus} as $\abs{a}:=\sqrt{b}$.
\edefin

An example of an operator that is not modular has to wait unit after the theorem below.

From
\beqn{
\AB{ax,ax}
~=~
\AB{x,bx}
~=~
\bAB{\abs{a}x,\abs{a}x}
}\eeqn
we see that $\ker a=\ker\abs{a}$. Clearly, $\ker\abs{a}\subset\ker b$, and from
\beqn{
\bAB{\abs{a}x,\abs{a}x}\,\bAB{\abs{a}x,\abs{a}x}
~=~
\AB{x,bx}\AB{bx,x}
~\le~
\AB{x,x}\,\norm{\AB{bx,bx}}
}\eeqn
we see that also $\ker b\subset\ker\abs{a}$, so $\ker b=\ker\abs{a}=\ker a$. Also, for every positive operator $b$ on a Hilbert module $E$, defining $E_b:=\ol{bE}$, we have:
\begin{itemize}
\item
$\ker b=E_b^\perp$ and $E_b\subset(\ker b)^\perp$. (Indeed, if $\AB{x,y}=0$ for all $y\in\ol{bE}$, then, in particular, $0=\AB{x,b^2x}=\AB{bx,bx}=0$ so $x\in\ker b$, and if $x\in\ker b$, then $\AB{x,y}=0$ for all $y\in\ol{bE}$. Summing up, $\ker b=E_b^\perp$. From this, we also get $E_b\subset E_b^{\perp\perp}=(\ker b)^\perp$.)

But, $E_b$ need not coincide with $(\ker b)^\perp$. (Indeed, suppose $E=\cB\ni\U$ and suppose $b\in\cB=\sB^a(E)$ is non-invertible but strictly positive. By definition, the latter means that $bx=0$ implies $x=0$, so $\ker b=\zero$ and $(\ker b)^\perp=E$. But, $E_b$ is the closed ideal generated by $b$, and since $b$ is not invertible, this ideal is non-trivial.)

\item
$E_b=E_{\sqrt{b}}$. (A special case of \cite[Proposition 3.7]{Lan95}.)

\item
$b$ is injective on $E_b$. (Indeed, let $x\in E_b=E_{\sqrt{b}}$, so there is a sequence $y_n\in E$ such that $x_n:=\sqrt{b}y_n\to x$. If $x\ne0$, then $0\ne\AB{x,x}=\lim_n\AB{y_n,\sqrt{b}x}$, so $\sqrt{b}x\ne0$, hence, $\AB{\sqrt{b}x,\sqrt{b}x}=\AB{x,bx}\ne0$, hence, $bx\ne0$.)
\end{itemize}
Now, if $a\in\sB^r(E,F)$ is modular, then we define $E_a:=E_b$ (for the unique $b$ illustrating that $a$ is modular). We get the following two results on conditions that assure a sort of polar decomposition more general than that known in literature.

\bthmn
Let $E$ and $F$ be Hilbert \nbd{\cB}modules and let $a\in\sB^r(E,F)$ be modular.

\begin{enumerate}
\item \label{1}
Considering $\abs{a}$ as an operator in $\sB^r(E,E_a)$, there is a unique isometry $v_a\in\sB^r(E_a,F)$ satisfying
\beqn{
v_a\abs{a}
~=~
a.
}\eeqn
Necessarily, $v_aE_a=\ol{aE}$.

\item \label{2}
There exists a partial isometry $v\in\sB^r(E,F)$ such that
\beqn{
vE
~=~
\ol{aE}
~~~~~~~~~
\text{and}
~~~~~~~~~
a
~=~
v\abs{a}
}\eeqn
if and only if $E_a$ is complemented in $E$. In that case, $v$ is unique.
\end{enumerate}
\ethmn

\proof
Clearly, the map $\abs{a}E\ni\abs{a}x\mapsto ax\in F$ is isometric, thus, extends to an isometry $v_a\colon E_a\rightarrow F$ satisfying everything (in particular, being uniquely determined by $v_a\abs{a}=a$) stated under \ref{1}.

Clearly, if there exists a partial isometry $v$ as in Part \ref{2}, then the restriction of $v$ to $E_a$ is bound to be $v_a$ from Part \ref{1}. Thanks to $vE=\ol{aE}$ we may (without changing, otherwise, the problem) pass to corestrictions of $a$ and $v$ to $F=\ol{aE}$ so that, now, $a$ has dense range and $v$ is surjective, hence, a coisometry. So, $E_a$ is the (recall, unique!) submodule such that $v$ restricts to unitary onto $\ol{aE}$, hence, complemented. So $E_a^{\perp\perp}=E_a$. Since $\ker v=E_a^\perp=\ker a$ it also follows that $v$ is bound to be $0$ on the remaining direct summand, so $v$ is unique.

If, conversely, $E_a$ is complemented, we may compose $v_a\colon E_a\rightarrow F$ from Part \ref{1} with the projection onto $E_a$ (considered as map $E\rightarrow E_a$) to get the partial isometry $v$ as stated.\qed

\bremn
In the usual treatment of polar decomposition, uniqueness is achieved by requiring $\ker v=\ker a$. We have $\ker a=E_a^\perp$, but we do not know in general if $(\ker a)^\perp=E_a^{\perp\perp}\supset E_a$ is actually equal to $E_a$. (See the example.) We wish to direct attention to how we replaced the condition $\ker v=\ker a$, which in the usual treatment achieves only uniqueness (but not here) while existence of $v$ is trivial, with the condition $vE=\ol{aE}$, which does not only achieve uniqueness but is also crucial for existence of $v$.
\eremn

\bobn
Suppose $a\in\sB^r(E,F)$ is modular and has range $aE$ dense in $F$. Then, similarly as in the proof of Part \ref{2}, not only $v$ (of which we do not yet know if it exists), but also $v_a$ (which by Part \ref{1} always exists) is surjective, hence, a unitary $v_a\colon E_a\rightarrow F$. Composing its adjoint $v_a^*\colon F\rightarrow E_a$ with the canonical embedding $E_a\rightarrow E$, we get an isometry $w\colon F\rightarrow E_a$, which is adjointable if and only if its range $E_a$ is complemented in $E$. If $w$ is adjointable, then $E_a=E_a^{\perp\perp}=(\ker a)^\perp$, so, if $(\ker a)^\perp=\zero$, then $E_a=\zero$, hence, $a=0$ (hence, $F=\zero$, because $aE=\zero$ is dense in $F$).
\eobn

A direct consequence of this observation is the following:

\bexn
With Kaad \cite{KaSk23}, on the standard Hilbert \nbd{\cB}module $E:=\cB^{\N_0}$ over the commutative, unital, and separable \nbd{C^*}algebra $\cB:=C\SB{0,1}$ we produced:
\begin{itemize}
\item
a surjective map $a\in\sB^r(E,\cB)$ (in fact, $a(b_0,0,\ldots)=b_0$);

\item
a submodule $F\subset E$ with $F^\perp=\zero$ and $\ker a\supset F$.

\end{itemize}
This means, in particular, that $(\ker a)^\perp\subset F^\perp=\zero$. If $a$ was modular, since $a$ is surjective, by the observation we could construct an isometry $w\colon\cB\rightarrow E$ onto $E_a$, and since $E_a=w\cB=(w\U)\cB=\xi\cB$ with the unit vector $\xi:=w\U$, we see that $E_a=\xi\xi^*E$ would be complemented in $E$. Since for complemented $E_a$ we have $E_a=E_a^{\perp\perp}=(\ker a)^\perp$, we would get $(\ker a)^\perp=E_a\ne\zero$. Consequently:

\noindent~\hfill\it{This map $a$ from \cite{KaSk23} is not modular.}\hfill~
\eexn

\bempn{Conclusion.~}
With the notion of not necessarily adjointable partial isometries from \cite{ShaSk23,Ske25} and the new notion of modular operators we did slightly enlarge applicability of the idea of polar decomposition for operators between Hilbert modules. Still, the applicability of polar decomposition for general \nbd{C^*}modules remains quite limited.
\eempn

\noindent
\bf{Acknowledgments:~}
I would like to thank the organizers of the ``\href{https://researchseminars.org/seminar/NYC-NCG}{Noncommutative geometry in NYC}'' online seminar, Alexander Katz and Igor Nikolaev, for giving me the occasion to talk about the paper \cite{Ske25}, of which this little note is a direct consequence. I thank Orr Shalit for several useful comments and hints to misprints.

\newpage

{
\setlength{\baselineskip}{2.5ex}


\begin{thebibliography}{BDH88}

\bibitem[BDH88]{BDH88}
M.~Baillet, Y.~Denizeau, and J.-F. Havet, \emph{{Indice d'une esperance
  conditionnelle}}, Compositio Math. \textbf{66} (1988), 199--236.

\bibitem[KS23]{KaSk23}
J.~Kaad and M.~Skeide, \emph{{Kernels of Hilbert module maps: A
  counterexample}}, J.\ Operator Theory \textbf{89} (2023), 343--348, (ar\-Xiv:
  2101.03030v1).

\bibitem[Lan95]{Lan95}
E.C. Lance, \emph{{Hilbert $C^*$--modules}}, Cambridge University Press, 1995.

\bibitem[Pas73]{Pas73}
W.L. Paschke, \emph{{Inner product modules over $B^*$--algebras}}, Trans.\
  Amer.\ Math.\ Soc. \textbf{182} (1973), 443--468.

\bibitem[Ske00]{Ske00b}
M.~Skeide, \emph{{Generalized matrix $C^*$--algebras and representations of
  Hilbert modules}}, Mathematical Proceedings of the Royal Irish Academy
  \textbf{100A} (2000), 11--38, (Cott\-bus, Rei\-he Mathe\-ma\-tik 1997/M-13).

\bibitem[Ske01]{Ske01}
\bysame, \emph{{Hilbert modules and applications in quantum probability}},
  Ha\-bi\-li\-ta\-tions\-schrift, Cottbus, 2001, available at:
  \\{\footnotesize\url{http://web.unimol.it/skeide/_MS/downloads/habil.pdf}}.

\bibitem[Ske06]{Ske06b}
\bysame, \emph{{Commutants of von Neumann correspondences and duality of
  Eilen\-berg-Watts theorems by Rieffel and by Blecher}}, Banach Center
  Publications \textbf{73} (2006), 391--408, (ar\-Xiv: math.OA/0502241).

\bibitem[Ske25]{Ske25}
\bysame, \emph{{Partial isometries between Hilbert modules and their
  compositions}}, Infin.\ Dimens.\ Anal.\ Quantum Probab.\ Relat.\ Top.
  \textbf{28} (2025), no.2, Paper 2440010 (16 pp), (ar\-Xiv: 2310.14755v2).

\bibitem[SS23]{ShaSk23}
O.M. Shalit and M.~Skeide, \emph{{CP-Semigroups and dilations, subproduct
  systems and superproduct systems: The multi-parameter case and beyond}},
  Dissertationes Math. \textbf{585} (2023), 1--233, (ar\-Xiv: 2003.05166v3).

\bibitem[WO93]{WeO93}
N.E. Wegge-Olsen, \emph{{$K$--Theory and $C^*$--algebras}}, Oxford University
  Press, 1993.

\end{thebibliography}

\newcommand{\Swap}[2]{#2#1}\newcommand{\Sort}[1]{}
\providecommand{\bysame}{\leavevmode\hbox to3em{\hrulefill}\thinspace}
\providecommand{\MR}{\relax\ifhmode\unskip\space\fi MR }
\providecommand{\MRhref}[2]{%
  \href{http://www.ams.org/mathscinet-getitem?mr=#1}{#2}
}
\providecommand{\href}[2]{#2}

}

\lf\noindent
Michael Skeide: \it{Dipartimento di Economia, Universit\`{a} degli Studi del Molise, Via de Sanctis, 86100 Campobasso, Italy},
E-mail: \href{mailto:skeide@unimol.it}{\tt{skeide@unimol.it}},\\
Homepage: \url{http://web.unimol.it/skeide/}


\end{document}